\documentclass[10pt,a4paper]{amsart}
\usepackage{verbatim}

\usepackage[T1]{fontenc}
\usepackage{graphicx}
\usepackage{enumerate}
\usepackage{amsmath,amsfonts,amssymb}
\usepackage{mathtools}
\usepackage{ulem}
\usepackage{color}
\usepackage{cite}
\definecolor{citation}{rgb}{0.2,0.58,0.2} 
\definecolor{formula}{rgb}{0.1,0.2,0.6}
\definecolor{url}{rgb}{0.3,0,0.5}
\usepackage{pgf,tikz}
\usepackage{mathrsfs}
\usepackage{fancyhdr}
\usepackage{dsfont}

\title[Measure data systems with  Orlicz growth]{Measure data systems with  Orlicz growth}

\author{Iwona Chlebicka}\address{Iwona Chlebicka \\
Institute of Applied Mathematics and Mechanics, University of Warsaw \\ ul. Banacha 2, 02-097 Warsaw, Poland\\  \texttt{e-mail: i.chlebicka@mimuw.edu.pl}} 

\author{Yeonghun Youn}\address{Yeonghun Youn\\ Department of Mathematics, Yeungnam University\\ Gyeongbuk 38541, Republic of Korea\\ \texttt{e-mail: yeonghunyoun@yu.ac.kr}}

\author{Anna Zatorska-Goldstein}\address{Anna Zatorska-Goldstein \\ Institute of Applied Mathematics and Mechanics, University of Warsaw \\ ul. Banacha 2, 02-097 Warsaw, Poland \\ \texttt{e-mail: azator@mimuw.edu.pl}}

\date{}
\begin{document}

\thanks{\vspace{1mm}}

\maketitle \sloppy

\thispagestyle{empty}

\belowdisplayskip=18pt plus 6pt minus 12pt \abovedisplayskip=18pt
plus 6pt minus 12pt
\parskip 4pt plus 1pt
\parindent 0pt

\newcommand{\ic}[1]{\textcolor{teal}{#1}}

\def\tens#1{\pmb{\mathsf{#1}}}
\newcommand{\barint}{
         \rule[.036in]{.12in}{.009in}\kern-.16in
          \displaystyle\int  } 
          
\newcommand{\dv}{{\rm div}}
\def\aI{\texttt{(a1)}}
\def\aII{\texttt{(a2)}}
\newcommand{\opA}{{\mathcal{ A}}}
\newcommand{\bopA}{{\bar{\opA}}}
\newcommand{\wt}{\widetilde}
\newcommand{\ve}{\varepsilon}
\newcommand{\vp}{\varphi}
\newcommand{\vt}{\vartheta}
\newcommand{\vk}{\varkappa}
\newcommand{\vr}{\varrho}
\newcommand{\pa}{\partial}
\newcommand{\cW}{{\mathcal{W}}}
\newcommand{\supp}{{\rm supp}}

\def\R{{\mathbb{R}}}
\def\N{{\mathbb{N}}}
\def\rp{{[0,\infty)}}
\def\r{{\mathbb{R}}}
\def\n{{\mathbb{N}}}
\def\l{{\mathbf{l}}}
\def\bu{{\bar{u}}}
\def\T{{\textsc{T}}}
\def\bg{{{g}}}
\def\bG{{{G}}}
\def\ba{{\bar{a}}}
\def\bv{{\bar{v}}}
\def\wtgamma{{\wt\gamma}}
\def\tew{{\tens{w}}}
\def\teeta{{\tens{\eta}}}
\def\texi{{\tens{\xi}}}
\def\teu{{\tens{u}}}
\def\teet{{\tens{\eta}}}
\def\tepsi{{\tens{\psi}}}
\def\teph{{\tens{\phi}}}
\def\tev{{\tens{v}}}
\def\telambda{{\tens{\lambda}}}
\def\calV{{\mathcal{V}}}
\def\tebv{{\bar{\tens{{v}}}}}
\def\tevp{{\tens{\vp}}}
\def\teF{{\tens{F}}}
\def\tef{{\tens{f}}}
\def\teg{{\tens{g}}}
\def\teelvr{{\tens{\ell_\vr}}}
\def\tewtu{{\tens{\wt u}}}
\def\tfA{{{\mathfrak{A}}}}

\def\teDu{{D\teu}}
\def\teDet{{D\teet}}
\def\teDph{{D\teph}}
\def\teDv{{D\tev}}
\def\teDw{{D\tew}}
\def\teDbv{{D\tebv}}
\def\teDwtu{{D\tewtu}}
\def\teDvp{D{\tevp}}
\def\teDelvr{{D\teelvr}}

\def\tedv{{\tens{\dv}}}
\def\temu{{\tens{\mu}}}
\def\tea{\tens{a}}

\def\tI{\text{I}}
\def\tII{\text{II}}
\def\tIII{\text{III}}
\def\bmu{{\bar{\mu}}}
\def\rn{{\mathbb{R}^{n}}}
\def\Rm{{\mathbb{R}^{m}}}
\def\Rn{{\mathbb{R}^{n}}}
\def\id{{\mathsf{Id}}}
\def\P{{\mathsf{P}}}
\def\Pxi{{\P_\texi}}
\def\Pj{{\mathsf{P_j}}}
\def\rnm{{\mathbb{R}^{n\times m}}}
\def\Rm{{\mathbb{R}^{m}}} 
\def\Mb{{\mathcal{M}(\Omega,\Rm)}} 

\newtheorem{coro}{\bf Corollary}[section]
\newtheorem{theo}[coro]{\bf Theorem} 
\newtheorem{lem}[coro]{\bf Lemma}
\newtheorem{rem}[coro]{\bf Remark} 
\newtheorem*{defi}{\bf Definition} 
\newtheorem{ex}[coro]{\bf Example} 
\newtheorem{fact}[coro]{\bf Fact} 
\newtheorem{prop}[coro]{\bf Proposition}

\newcommand{\apdu}{{\rm ap} \, D \teu}
\newcommand{\temf}{f} 
\newcommand{\data}{\textit{\texttt{data}}}


\parindent 1em

\begin{abstract}
We study the existence of very weak solutions to a system
\[\begin{cases}-\tedv \opA(x,\teDu)=\tens{\mu}\quad\text{in }\ \Omega,\\
    \tens{u}=0\quad\text{on }\ \partial\Omega\end{cases}
\]
 with a datum $\temu$ being a vector-valued bounded Radon measure and $\opA:\Omega\times\rnm\to\rnm$ having measurable dependence on the spacial variable and Orlicz growth with respect to the second variable.  We are {\em not} restricted to the superquadratic case.
 
 For the solutions and their gradients we provide regularity estimates in the generalized Marcinkiewicz scale. In addition, we show a precise sufficient condition for the solution to be a~Sobolev function.
 
\end{abstract}
\section{Introduction}

We concentrate on the problem
\[\begin{cases}-\tedv \opA(x,\teDu)=\tens{\mu}\quad\text{in }\ \Omega,\\
    \tens{u}=0\quad\text{on }\ \partial\Omega\end{cases}
\]
involving nonlinear operator $\opA:\Omega\times\rnm\to\rnm$ with measurable dependence on $x$ and having Orlicz growth with respect to the second variable. Our main model is
 \begin{equation}
\label{intro:eq:main} -{ \tedv}\left(a(x)\frac{G(|\teDu|)}{|\teDu|^2}  \teDu \right)=\temu\quad\text{in}\quad \Omega,
\end{equation}
where $ a\in L^\infty(\Omega)$ is a scalar function separated from zero, $\temu$ is a vector-valued bounded Radon measure, and $G\in\Delta_2\cap\nabla_2$. For details  see {\rm Assumption {\bf (A-ex-vect)}} below.  We cover the case of $\opA$ being $p$-Laplacian when $G_p(s)=s^{p},$ for every $p>1$, but also operators governed by the Zygmund-type functions $G_{p,\alpha}(s)=s^p\log^\alpha(1+s)$, $p>1,\,\alpha\in\R$, as well as their multiplications and compositions with various parameters. {We study the existence of very weak solutions using a relevant notion of solution to the system. Moreover we  provide regularity estimates in the generalized Marcinkiewicz scale for the solutions and their gradients.} We show a~sharp growth condition on $G$ sufficient for the solution to belong to $W^{1,1}_{loc}(\Omega,\Rm)$.\newline

Since weak solutions do not have to exist for arbitrary measure datum, a~weaker notion of solution should be employed.  {Distributional solutions can be non-unique and almost as bad as any function from a Sobolev space, e.g. unbounded on every open subset}~\cite{Serrin-pat}. Various kinds of very weak solutions to measure data nonlinear equations of power growth has been introduced since the fundamental contributions~\cite{BG,bbggpv}. Following their ideas, we consider so-called approximable solutions and SOLA (Solutions Obtained as a Limit of Approximation). 
Recently the theory has been developed in Orlicz and generalized Orlicz setting e.g. by \cite{ACCZG,CGZG,IC-measure-data,CiMa}. For pioneering work on the  equations with Orlicz growth we refer to~\cite{G,lieb,T}, while for regularity to their measure data counterparts see~\cite{Baroni-Riesz,Y-new,IC-gradest,IC-lower,CGZG-Wolff}. The cornerstone of studies on quasilinear elliptic systems in divergence form has been laid in~\cite{Ural,Uhl}. Existence of weak solutions to systems with Orlicz growth was proven in~\cite{DS}. The studies on elliptic {systems} involving measure data started with~\cite{FR,Le,Rak} and were followed by preeminent works~\cite{DHM,DHM-n,DHM-gen}, s
ee also~\cite{L1,LP1,LP2,LP3,LRS1,LRS2,zhou}.  Simultaneously, the theory of regularity to such problems, as well as to related problems in the calculus of variations was being developed in many contexts. For vectorial $p$-Laplace problems see e.g. \cite{Byun5,CiSch,KuMi2016p,KuMi2018}, for doubling Orlicz growth problems see e.g. \cite{BCDM,CiMa-ARMA2014,DiEt,DiLeStVe,DF1}, while for problems under other growth conditions see e.g.  ~\cite{bemi,CDFK,DF-vec,DFMin,DiMa,Marc2006,MarcPapi}. For broader overview see~\cite{IC-pocket,KuMi2016,Marc2020,Marc2021,MiRa}. \newline

Let us give the details about the measure data problem we study.\newline

\noindent{\em Essential notation}. By {`$\cdot$'} we denote the {scalar product of two vectors, i.e. for ${\bf{\xi}}=(\xi_1,\dots,\xi_n)\in \R^n$ and 
	${\bf{\eta}}= (\eta_1,\dots,\eta_n)\in \R^n$ we have ${\bf{\xi}}\cdot{\bf{\eta}} = \sum_{i=1}^n \xi_i \eta_i$};\\ by {`$:$'} -- {the Frobenius product of the second-order tensors, i.e. for $\tens{\xi}=[\xi_{j}^\alpha]_{j=1,\dots,n,\, \alpha=1,\dots,m}$ and $\tens{\eta}=[\eta_{j}^\alpha]_{j=1,\dots,n,\, \alpha=1,\dots,m}$ we have
	\[\tens{\xi}: \tens{\eta} =\sum_{\alpha=1}^m \sum_{j=1}^n \xi_{j}^\alpha \eta_{j}^\alpha.\]}
	By `{$\otimes$}' we denote {the tensor product of two vectors, i.e for ${{\xi}}=(\xi_1,\dots,\xi_n)\in \rn$ and 
	${{\eta}}= (\eta_1,\dots,\eta_n)\in \rn$, we have $\xi\otimes\eta:=[\xi_i\eta_j]_{i=1,\dots,n,\,j=1,\dots,n},$ that is \[{\xi}\otimes{\eta}:=\begin{pmatrix}
  \xi_{1}\eta_{1} & \xi_{1}\eta_{2} & \cdots & \xi_{1}\eta_{n} \\
  \xi_{2}\eta_{1} & \xi_{2}\eta_{2} & \cdots & \xi_{2}\eta_{n} \\
  \vdots  & \vdots  & \ddots & \vdots  \\
  \xi_{n}\eta_{1} & \xi_{n}\eta_{2} & \cdots & \xi_{n}\eta_{n} 
 \end{pmatrix}\in \R^{n\times n}.\]}
 
\noindent{\em Problem}.  Given a bounded open set $\Omega\subset\rn$, $n\geq 2$, we investigate solutions $\teu:\Omega\to\Rm$ to the problem \begin{equation}
    \label{eq:mu}\begin{cases}-\dv \opA(x,\teDu)=\temu\quad\text{in }\ \Omega,\\
    \teu=0\quad\text{on }\ \partial\Omega\end{cases}
\end{equation}
 with a datum $\temu$ being a vector-valued bounded Radon measure and $\opA:\Omega\times\rnm\to\rnm$ which satisfies growth and coercivity conditions expressed by the means of an $N$-function $G$. We admit $G$ being a Young function such that $g=G'$ satisfies the following condition
\begin{equation}\label{iG-sG} 
0 < i_g=\inf_{t>0}\frac{tg'(t)}{g(t)}\leq \sup_{t>0}\frac{tg'(t)}{g(t)}=s_g<\infty.\end{equation} We consider a problem~\eqref{eq:mu} under the following regime.\newline

\noindent{\em Assumption {\bf (A-ex-vect)}}. Assume that $\opA: {\Omega \times} \rnm\to\rnm$ is {a Carath\'eodory's function which is} strictly monotone in the sense that\begin{equation}
    \label{ass:str-mon}\big(\opA(x,\xi)-\opA(x,\eta)\big):(\xi-\eta)>0\qquad\text{whenever }\xi\neq\eta.
\end{equation}
Moreover, suppose that $\opA(x,0)=0$ and $\opA$ satisfies the following growth and coercivity conditions\begin{equation}
    \label{ass:gr'n'coer} \begin{cases}\opA(x,\xi):\xi&\geq c_1^\opA G(|\xi|),\\
    |\opA(x,\xi)|&\leq c_2^\opA\left(g(|\xi|)+b(x)\right)\end{cases}
\end{equation}
for a.a. $x\in\Omega$, all $\xi\in\rnm,$ and some $b\in L^{\wt G}(\Omega).$ Furthermore, we require $\opA$ to  satisfy\begin{equation}
    \label{ass:struct} \opA(x,\xi):\big((\id-w\otimes w)\xi\big)\geq 0
\end{equation}
for a.a. $x\in\Omega$, all $\xi\in\rnm,$ and every vector $w\in\Rm$ with $|w|\leq 1$. 
\newline

{Condition~\eqref{ass:struct} comes from~\cite{DHM-gen}. It is satisfied by any $\opA(x,\xi)=\alpha(x,\xi)\xi$ where $\alpha$ is a real valued non-negative function, which justifies to call system involving it as  {\em quasidiagonal}. In particular, it is satisfied in the case of $p$-Laplacian or~\eqref{intro:eq:main}. By now~\eqref{ass:struct} is a typical assumption in the analysis of elliptic systems. We use it to control the energy of approximate solutions, e.g. in order to get uniform a priori estimates. }\newline

We shall deal with two kinds of very weak solutions obtained in the approximating procedure motivated by~\cite{BG} that differ in the expected regularity and the required convergence of approximate problems.

\begin{defi}[Approximable solution] \label{def:as:mu}
 A vector valued map $\teu\in  {\Upsilon}^{1,G}(\Omega,\Rm)$  is called an approximable solution to~\eqref{eq:mu} under the regime of {\rm Assumption {\bf (A-ex-vect)}}, if there exists a sequence $(\teu_{h})\subset W^{1,G }(\Omega,\Rm)$ of local energy solutions to the systems 
\begin{equation}
    \label{eq:app-prob}-\tedv\opA(x,D\teu_{h})={\temu_h}
\end{equation} such that 
\[ \teu_h\to \teu \quad\text{and}\quad \opA(x, \teDu_h )\rightarrow \opA(x,\teDu) \quad \hbox{a.e. in  $\Omega$}\]
 and $(\temu_h)\subset L^\infty (\Omega,\Rm)$ is a sequence
of maps that converges to $\temu$ weakly in the sense of measures and satisfies
\begin{equation}
    \label{conv-of-meas}
\limsup_h |\temu_h |(B) \leq |\temu|(B)\qquad\text{
for every ball $B\subset\Omega$.}
\end{equation}
 {For the definition of a space $\Upsilon^{1,G}(\Omega,\Rm)$ see Section 2.3}
\end{defi}

\begin{defi}[SOLA]\label{def:appr-sol}
 A vector valued map $\teu\in W_0^{1,1}(\Omega,\Rm)$ such that $\int_\Omega g(|\teDu|)\,dx<\infty$ is called a SOLA to~\eqref{eq:mu} under the regime of {\rm Assumption {\bf (A-ex-vect)}}, if there exists a sequence $(\teu_{h})\subset W^{1,G }(\Omega,\Rm)$ of local energy solutions to the systems~\eqref{eq:app-prob} such that \[\text{$\teu_{h}\to \teu\quad$ locally in $\quad W^{1,1}(\Omega,\Rm)$}\] and $(\temu_h)\subset L^\infty (\Omega,\Rm)$ is a sequence
of maps that converges to $\temu$ weakly in the sense of measures and satisfies~\eqref{conv-of-meas}.
\end{defi}
{Observe that the above approximation property immediately implies that a SOLA $\teu$ is a distributional solution to~\eqref{eq:mu}, that is,
\[\int_\Omega\opA(x,\teDu):\teDvp\,dx=\int\tevp \,d\temu\qquad\text{for every }\tevp\in  C_0^\infty (\Omega,\Rm).\]}

For $G_n$ being a Sobolev conjugate to $G$ given by~\eqref{GN} let us define\begin{equation}
    \label{thetatheta}\vt_n(t):=\frac{G_n(t^\frac{1}{n'})}{t}\qquad\text{ and }\qquad \theta_n(t):=\frac{t}{G_n^{-1}(t)^{n'}}\,.
\end{equation}
Let us state a result on existence and Marcinkiewicz regularity of approximable solutions to problems with general structure. We distinguish two cases related to different Sobolev-type embeddings for slowly and fast growing $G$ described in Section~\ref{ssec:emb}, see~\eqref{intG}.
\begin{theo} \label{theo:exist}
If a vector field $\opA$ satisfies {\rm Assumption {\bf (A-ex-vect)}} and $\temu\in\Mb$, then there exists an approximable solution $\teu\in {\Upsilon}^{1,G}_0(\Omega,\Rm)$ to~\eqref{eq:mu} for which \begin{equation}
    \int_\Omega g(|\teDu|)\leq c|\temu|(\Omega)\qquad\text{with some }\ c=c(\data)>0.
\end{equation}
Moreover\begin{equation}
    \label{distr}\int_\Omega\opA(x,\teDu):\teDvp\,dx=\int_\Omega\tevp\,d\temu\qquad\text{for any }\ \tevp\in W^{1,\infty}_0(\Omega,\Rm).
\end{equation}
Furthermore the following assertions hold true.
\begin{itemize}
    \item[(i)]If $G$ grows so fast that \eqref{intG}$_2$ holds for  $\sigma = n$,  approximable solution $\teu$ is a~weak solution and $\teu\in W_0^{1,G}(\Omega,\Rm)\subset L^\infty(\Omega)$.
    \item[(ii)] If $G$ grows so slowly that \eqref{intG}$_1$ holds, we have
\begin{equation}
    \label{Marc-reg}
    |\teu|\in L^{\vt_n(\cdot),\infty}(\Omega)\qquad\text{and}\qquad 
    |\teDu|\in L^{\theta_n(\cdot),\infty}(\Omega).
\end{equation}
\end{itemize}
\end{theo}
In order to prove that approximable solutions are Sobolev functions  
we need to impose a growth condition on $G$. Let us take $H_n$ defined in~\eqref{GN} with $\sigma=n$ and denote
\begin{equation*}
     \Psi_{n}(t) = \frac{G(t)}{H_{n}(t)^{n'}}.
\end{equation*}
Note that $\Psi_n = \theta_n \circ G$.

\begin{theo}\label{theo:SOLA} Suppose a vector field $\opA$ satisfies {\rm Assumption {\bf (A-ex-vect)}}, $\temu\in\Mb$ and $G$ grows fast enough to satisfy
    \begin{equation}\label{SOLA-Psi}
    \int^{\infty} \frac{dt}{\Psi_{n}(t)} < \infty.
    \end{equation}
Then each approximable solution $\teu$ to~\eqref{eq:mu}  satisfies $\teu\in  W^{1,1}_{0}(\Omega,\Rm)$ and $\int_\Omega g(|\teDu|)\,dx<\infty$, hence it is a SOLA.
\end{theo}
\begin{coro}
Suppose a vector field $\opA$ satisfies {\rm Assumption {\bf (A-ex-vect)}}, $\temu\in\Mb$ and $G$ grows fast enough to satisfy~\eqref{SOLA-Psi}, then there exists a SOLA to~\eqref{eq:mu} and it shares Marcinkiewicz-type regularity of {\rm (i)-(ii)} of Theorem~\ref{theo:exist}.
\end{coro}
\begin{rem}[Sharpness]\rm
Condition \eqref{SOLA-Psi} in the case of $G(t)=t^p$, $p>1$ reads $p>2-1/n$, cf.~\cite{BG}. This bound cannot be extended even in the case of a single equation. Indeed, it is enough to consider the fundamental solution, i.e. $\bar u=|x|^\frac{p-n}{p-1}$ being distributional solution to $-\Delta_p \bar u=\delta_0$ with $1<p<n$ and $\delta_0$ being a Dirac delta. It  has locally integrable gradient if and only if  $p>2-1/n$.  For the Orlicz counterpart of this fact see~\cite[Corollary~2.4]{CGZG-Wolff}.
\end{rem}
\begin{rem}[Regularity]\rm Known method giving precise information on very solutions to measure data systems of a form~\eqref{eq:mu} and power growth are potential estimates, see~\cite{KuMi2016,KuMi2018}. For an Orlicz growing operator of quasidiagonal structure and superquadratic growth~\cite{CYZG} provides an estimate on the solution by the means of potential of a Wolff type for and its regularity consequences including precise continuity and H\"older continuity criteria. 
\end{rem}

\noindent{\em Our methods}. The considered definitions of very weak solutions and general outline of our proof of their existence base on the approximate procedure coming from pioneering papers~\cite{BG,bbggpv}. We follow the techniques introduced there and in the already classical reasoning for measure data systems~\cite{DHM}. However, general growth of our operator require more delicate approach. Let us recall that we do not have to impose any growth conditions from below for the existence of approximable solutions. In order to obtain certain a priori estimates we employ ideas of a very recent paper~\cite{BCDM}. Getting regularity in the Marcinkiewicz-type scale is proven by applying optimal Orlicz--Sobolev embedding as in the proof in~\cite{CiMa} for scalar Orlicz problems. In order to find a sufficient condition for the solution to be a~Sobolev function we use the method of~\cite{Y-new}.  

\section{Preliminaries}\label{sec:prelim}

\subsection{Notation}
In the following, we shall adopt the customary convention of denoting by $c$ a constant that may vary from line to line. Sometimes to skip rewriting a constant, we use $\lesssim$. By $a\approx b$, we mean $a\lesssim b$ and $b\lesssim a$. To stress the dependence of the intrinsic constants on the parameters of the problem, we write $\lesssim_\data$ or $\approx_\data.$ 
We make use of the truncation operators at the level $k$. One-dimensional version is denoted by $\T_k:\R\to\R$ and given by $\T_k(f):=\min\left\{1,\frac{k}{|f|}\right\}f$, while multi-dimesional one $T_k:\Rm\to\Rm$ is defined as follows 
\begin{equation}\label{Tk}T_k(\texi):=\min\left\{1,\frac{k}{|\texi|}\right\}\texi.
\end{equation}
Then, of course, $DT_k:\Rm\to\Rm$ is given by\begin{equation}
    \label{DTk}DT_k(\texi)=\begin{cases}
    \id&\text{if }\ |\texi|\leq k,\\
    \frac{k}{|\texi|}\left(\id-\frac{\texi \otimes \texi}{|\texi|^2}\right)&\text{if }\ |\texi|> k.
    \end{cases}
\end{equation} 
 {We shall also consider $\Theta_t \in C^1_0(\Rm,\Rm)$, which satisfies \begin{equation}
    \Theta_t =\begin{cases}{\id} &\text{ on $B(0,t)$},\\
    \text{with $D\Theta_t$ bounded}& \text{ on  $B(0,2t)\setminus B(0,t)$},\\ 
    0 &\text{ outside $B(0,2t)$}\end{cases} \label{Theta-t}
\end{equation}
and a set\begin{equation}
    \label{M-alpha}
M_{\alpha}= \{ x \in \Omega :\ |\teu(x)| <\alpha \}.
\end{equation}}

To describe ellipticity of a vector field $\opA$, we make use of a function $\calV:\rnm\to\rnm$
\begin{equation}
    \label{V-def}
\calV({\xi})= \left(\frac{g(|{\xi}|)}{|{\xi}|}\right)^{1/2}\xi . 
\end{equation}
For the case of referring to the dependence of some quantities on the parameters of the problem, we collect them as  \[\data=\data(c_1^\opA,c_2^\opA,i_g,s_g,\|b\|_{L^{\wt G}},n,m).\]
Note that $c_1^\opA,c_2^\opA,b$ are given by~\eqref{ass:gr'n'coer}, $i_g,s_g$ describe the growth of $g$ -- see~\eqref{iG-sG}, while $n,m$ are dimensions of the system we consider.

\subsection{Basic definitions} References for this section {are~\cite{rao-ren,KR}}.

We say that a function $G: [0, \infty) \to [0, \infty]$ {is a Young function} if it is convex, vanishes  at $0$, and is neither identically equal to 0, nor to infinity. A Young function $G$ which is finite-valued,
vanishes only at $0$ and satisfies the additional growth conditions
\begin{equation*}\lim _{t \to
0}\frac{G (t)}{t}=0 \qquad \hbox{and} \qquad \lim _{t \to \infty
}\frac{G (t)}{t}=\infty 
\end{equation*} 
is called an $N$-function.  The  complementary~function $\wt{G}$  (called also the Young conjugate, or the Legendre transform) to a nondecreasing function $G:\rp\to\rp$  is given by the following formula
\[\wt{G}(s):=\sup_{t>0}(s\cdot t-G(t)).\]
If $G$ is a Young function, so is $\wt{G}$. 
If $G$ is an $N$-function, so is $\wt{G}$.

Having Young functions $G,\wt G$, we are equipped with Young's inequality reading as follows\begin{equation}
\label{in:Young} ts\leq G(t)+\wt{G}(s)\quad\text{for all }\ s,t\geq 0.
\end{equation}

 We say that a function $G:\rp\to\rp$ satisfies $\Delta_2$-condition if there exist $c_{\Delta_2},t_0>0$ such that $G(2t)\leq c_{\Delta_2}G(t)$ for $t>t_0.$ 
 We say that $G$ satisfy $\nabla_2$-condition if $\wt{G}\in\Delta_2.$ 
Note that it is possible that $G$ satisfies only one of the conditions $\Delta_2/\nabla_2$. For instance, for $G(t) =( (1+|t|)\log(1+|t|)-|t|)\in\Delta_2$, its complementary function is  $\widetilde{G}(s)= (\exp(|s|)-|s|-1 )\not\in\Delta_2$. See~\cite[Section~2.3, Theorem~3]{rao-ren} for equivalence of various definitions of these conditions and~\cite{CGZG,DFL} for illustrating the subtleties. In particular, $G\in\Delta_2\cap\nabla_2$ if and only if $1<i_G\leq s_G<\infty$. This assumption implies a comparison with power-type functions i.e.
$\frac{G(t)}{t^{i_G}}$ 
 is non-decreasing 
  and 
  $\frac{G(t)}{t^{s_G}}$ 
  is non-increasing, and it is stronger than being sandwiched between power functions.
 
 \begin{lem}\label{lem:equivalences}
 If {an $N$-function} $G\in\Delta_2\cap\nabla_2$, then 
 $g(t)t\approx  G(t)$ 
  and 
 $\wt G(g(t))\approx G(t)$ 
 with the constants depending only on the growth indexes of $G$, that is $i_G$ and $s_G$.  Moreover,  $g^{-1}(2t)\leq c g^{-1}(t)$ with $c=c(i_G,s_G).$
 \end{lem}

Due to Lemma~\ref{lem:equivalences} and~\cite[Lemmas~3 and~21]{DiEt}, we have the following relations.
\begin{lem}\label{lem:DiEt-mon} If $G$ is an $N$-function of class $C^2((0,\infty))\cap C([0,\infty))$, $g$ satisfies~\eqref{iG-sG}, and $\opA$ satisfies~\eqref{ass:gr'n'coer}, then for every $\xi,\eta\in\rnm$ it holds 
\begin{equation}
    \label{opA:strict-monotonicity}
\left(\opA(x,\xi)-\opA(x,\eta)\right):{(\xi-\eta)}\gtrsim_{\data}\, \frac{g(|\xi|+|\eta|)}{|\xi|+|\eta|}|\xi-\eta|^2\approx_{\data} \left| \calV(\xi)-\calV(\xi) \right|^2,
\end{equation}
and\begin{equation}
    \label{relation:g-V} g(|\xi|+|\eta|)|\xi-\eta|\approx_{\data} G^\frac{1}{2}(|\xi|+|\eta|) |\calV(\xi)-\calV(\eta)|.
\end{equation}
\end{lem}
 \subsection{ {Notions of gradients and function} spaces}
 Basic reference  for this section is~\cite{adams-fournier}, where the theory of Orlicz spaces is presented for scalar functions. The proofs for functions with values in $\Rm$ can be obtained by obvious modifications.
 
We study the solutions to PDEs in the Orlicz--Sobolev spaces equipped with a~modular function $G\in C^1 {((0,\infty))}$ -- a strictly increasing and convex function   such that $G(0)=0$ and satisfying~\eqref{iG-sG}. Let us define a modular \begin{equation}
    \label{modular}
\vr_{G,U}(\texi)=\int_U G(|\texi|)\,dx.
\end{equation}
 
For any bounded $\Omega\subset\rn$, by Orlicz space ${L}^G(\Omega,\Rm)$, we understand the space of measurable functions endowed with the Luxemburg norm 
\[||\tef||_{L^G(\Omega, \Rm)}=\inf\Big\{\lambda>0:\ \ \vr_{G,\Omega}\left( \tfrac{1}{\lambda} |\tef|\right)\leq 1\Big\}.\]
The counterpart of the H\"older inequality in Orlicz spaces reads \begin{equation}
\label{in:Hold} \|\tef\teg\|_{L^1(\Omega,\Rm)}\leq 2\|\tef\|_{L^G(\Omega,\Rm)}\|\teg\|_{L^{\wt{G}}(\Omega,\Rm)}
\end{equation}
for all $\tef\in L^G(\Omega,\Rm)$ and $ \teg\in L^{\wt{G}}(\Omega,\Rm)$.

 We define the Orlicz--Sobolev space  $W^{1,G}(\Omega, \Rm)$  as follows
\begin{equation*} 
W^{1,G}(\Omega,\Rm)=\Big\{\tef\in W^{1,1}_{loc}(\Omega,\Rm):\ \ |\tef|,|D{\tef}|\in L^G(\Omega,\Rm)\Big\},
\end{equation*}where the gradient is understood in the distributional sense, endowed with the norm
\[
\|\tef\|_{W^{1,G}(\Omega,\Rm)}=\inf\Big\{\lambda>0 :\ \    \vr_{G,\Omega}\left( \tfrac{1}{\lambda} |\tef|\right)+ \vr_{G,\Omega}\left( \tfrac{1}{\lambda} |D{\tef}|\right)\leq 1\Big\} 
\]
and  by $W_0^{1,G}(\Omega,\Rm)$ we denote the closure of $C_c^\infty(\Omega,\Rm)$ under the above norm. 


 A measurable function $\teu$ is said to be approximately differentiable at $x\in\Omega$ if there exists a $n\times m$ matrix $F_x$ such that for every $\ve>0$ it holds
\begin{equation}
    \label{app-grad-def}
\lim_{r\to 0}
\frac{\left|\left\{y\in B(x,r):\ |\teu(y)-\teu(x)-F_x(y-x)|>\ve r\right\}\right|}{r^n}=0\,.
\end{equation}
We write $F_x={\rm ap} \teDu(x)$. For more information see \cite{EG}.

We write that $\teu\in{\Upsilon}_0^{1,G}(\Omega,\Rm)$ if for every  {$\Theta \in C^1_0(\Rm,\Rm)$ we have $\Theta(\teu)\in W^{1,G}_0(\Omega,\Rm)$.}  
{\begin{lem}
If $\teu\in {\Upsilon}_0^{1,G}(\Omega,\Rm)$, then there exists a unique measurable function
$\tens{Z_u} : \Omega \to \rnm$ such that for $\Theta_t \in C^1_0(\Rm,\Rm)$  which satisfies \eqref{Theta-t}, with arbitrary $t>0$, it holds 
\begin{equation}\label{gengrad'} 
D \Theta_t(\teu) = \tens{Z_u} \,\, \hbox{
a.e. in $\{|\teu|<t\}$ 
for every $t > 0$}
\end{equation} 
and \begin{equation}\label{gengrad-lem'} 
D T_t(\teu) = \tens{Z_u} \,\, \hbox{
a.e. in $\{|\teu|<t\}$ 
for every $t > 0$.}
\end{equation}
\end{lem}
\begin{proof}
Observe that arguments of \cite[Lemma~2.1]{bbggpv} work also in the vectorial case apart from the fact that the gradient of vectorial truncation $T_t(\teu)$ does not vanish outside the set $\{|\teu|<t\}$ 
for every $t > 0$. Note however that for $\ve>0$ it holds
\[D\Theta_t(\Theta_{t+\ve}(\teu))=D\Theta_t(\teu)=\tens{Z_u} \,\, \hbox{
a.e. in $\{|\teu|<t\}$ 
for every $t > 0$.}\]
Note that \eqref{gengrad-lem'} follows from the same arguments.
\end{proof}
 }

 One has that $\teu\in W_0^{1,G}(\Omega,\Rm)$ if and only
if $\teu\in {\Upsilon}_0^{1,G}(\Omega,\Rm)$
and
$\tens{Z_u}\in L^G(\Omega, \rnm)$. In the latter case, $\tens{Z_u} = \teDu$ a.e. in $\Omega$. With an abuse of~notation, for every $\teu\in
 {\Upsilon}_0^{1, G}(\Omega,\Rm)$ we denote $\tens{Z_u}$ simply by $\teDu$ throughout.

 We will use that if $G$ is an $N$-function, a bounded sequence in $L^G(\Omega,\rnm)$ is uniformly integrable in $L^1(\Omega;\rn)$, which can be interpreted by the classical Dunford-Pettis' theorem.
\begin{theo}[Dunford-Pettis]\label{theo:dunf-pet}
A family $\{U_k\}_{k}$ of measurable functions  is uniformly integrable in $L^1(\Omega;\rn)$ if and only if it is relatively compact in the weak topology.
\end{theo}

\subsection{Embeddings}\label{ssec:emb}
For Sobolev--Orlicz spaces expected embedding theorems hold true. We apply \cite[Theorem~4.1]{BCDM}. For this we consider $\sigma\geq n$ and we distinguish two possible behaviours of $G$
\begin{equation}
\label{intG}
\int^\infty\left(\frac{t}{G(t)}\right)^\frac{1}{\sigma-1}dt=\infty \qquad\text{and}\qquad
\int^\infty\left(\frac{t}{G(t)}\right)^\frac{1}{\sigma-1}dt<\infty,
\end{equation}
which roughly speaking describe slow and fast growth of $G$ in infinity, respectively. This result is proven under the restriction 
\begin{equation}\label{int0G}\int_0\left(\frac{t}{G(t)}\right)^\frac{1}{\sigma-1}dt<\infty, 
\end{equation} 
concerning the growth of $G$ in the origin. Nonetheless, the properties of $L^G$ depend on the behaviour of $G(s)$ for large values of $s$ and~\eqref{int0G} can be easily by-passed in application by considering a function that generates the same function space \begin{equation}
    \label{modification-at-0}
{G}^0(t)=tG(1)\mathds{1}_{[0,1]}(t)+G(t)\mathds{1}_{(1,\infty)}(t).
\end{equation} Indeed, $L^G(\Omega)=L^{G^0}(\Omega).$ If $G$ is a Young function satisfying \eqref{int0G} and \eqref{intG}$_1 $, consider 
\begin{equation}\label{GN}
H_\sigma(s)=\left(\int_0^s\left(\frac{t}{G(t)}\right)^\frac{1}{\sigma-1}dt\right)^\frac{\sigma-1}{\sigma} \qquad \text{and}\qquad
G_\sigma(t)=G(H_\sigma^{-1}(t)).
\end{equation}
Then due to~\cite[Theorem~4.1]{BCDM}  there exists a constant $c_s=c_s(n,\sigma,|\Omega|)$, such that for every $u\in W_0^{1,G}(\Omega)$ it holds that \begin{equation}
    \label{Sob-slow}\int_\Omega G_\sigma\left(\frac{|u|}{c_s\big(\int_\Omega G(|Du|)dx\big)^\frac{1}{\sigma}}\right)dx\leq \int_\Omega G(|Du|)dx\quad\text{for every $u\in W_0^{1,G}(\Omega)$.}
\end{equation}

If $G$ is a Young function satisfying \eqref{intG}$_2 $ for $n=\sigma$, then by~\cite{AC} we have $ W_0^{1,G}(\Omega)\subset {L^\infty(\Omega)}$. Moreover, when we denote
\begin{equation}
    \label{Fn}F_n(t)=t^{n'}\int_t^\infty \frac{\wt G(s)}{s^{1+n'}}\,ds\qquad \text{and}\qquad
J_n(s)=\frac{s}{F_n^{-1}(s)},\end{equation} we may write that there exists a constant $c=c(n,|\Omega|)$, such that
\begin{equation}
    \label{Sob-fast}\|u\|_{L^\infty(\Omega)}\leq c J_n\left(\int_\Omega G(|Du|)\,dx\right) \quad\text{for every $u\in W_0^{1,G}(\Omega)$.}
\end{equation}

\subsection {The operator and definition of solutions} We study the problem~\eqref{eq:mu} with $\opA$ as in {\rm Assumption {\bf (A-ex-vect)}}. We notice that under such regime  the operator   $\mathsf{A}_{G}$ acting as
\begin{flalign*}
\langle\mathsf{A}_{G}\teu,\teph\rangle:=\int_{\Omega}\opA(x,\teDu): \teDph \, dx\quad \text{for}\quad \teph\in C^{\infty}_{0}(\Omega,\Rm)
\end{flalign*}
is well-defined on a reflexive and separable Banach space $W_0^{1,G}(\Omega,\Rm)$ and $\mathsf{A}_{G}(W_0^{1,G}(\Omega,\Rm))\subset (W_0^{1,G}(\Omega,\Rm))'$. Indeed, when $\teu\in W_0^{1,G}(\Omega,\Rm)$ and $\teph\in C_c^\infty(\Omega,\Rm)$, structure condition~\eqref{ass:gr'n'coer},  H\"older's inequality~\eqref{in:Hold}  justify that
\begin{flalign*}
|{\langle \mathsf{A}_{G}\teu,\teph \rangle}|\le &\,c\,\int_{\Omega}g(|{\teDu}|){|\teDph|} \, dx \le c\left \| g(|{\teDu}|) \right \|_{L^{\wt G(\cdot)} }\|{|\teDph|}\|_{L^{G}}.\nonumber 
\end{flalign*}
Since for  $\teu\in W_0^{1,G}(\Omega,\Rm)$  we have $\left \| g(|{\teDu}|) \right \|_{L^{\wt G(\cdot)} }<c$, we can conclude that
\begin{flalign*}
|{\langle \mathsf{A}_{G}\teu,\teph \rangle}| &\leq\, c\|{|\teDph|}\|_{L^{G}}\le c\|{\teph}\|_{W^{1,G}}.
\end{flalign*}

A function $\teu\in W^{1,G}_{loc}(\Omega,\Rm)$ is called {a} {\em weak solution} to~\eqref{eq:mu}, if\begin{equation}
    \label{eq:main-mu-weak}\int_\Omega \opA(x,\teDu): \teDvp\,dx=\int_\Omega\tevp\,d\temu(x)\quad\text{for every }\ \tevp\in W^{1,G}_{0}(\Omega,\Rm) \,.
\end{equation}
Recall that $W^{1,G}_{0}(\Omega,\Rm)$ is separable and by its very definition ${C_c^\infty}(\Omega,\Rm)$ is dense there.
 
\begin{rem}[Existence and uniqueness of weak solutions] \rm For $\temu\in (W^{1,G}_{0}(\Omega,\Rm))',$ due to the strict monotonicity of the operator, there exists a unique weak solution to~\eqref{eq:mu}, see~\cite[Section~3.1]{KiSt}.
\label{rem:weak-sol}
\end{rem}

\section{Existence and Marcinkiewicz regularity to vectorial measure data problems}\label{sec:ex-vec}

In this section we prove our main accomplishments, i.e. Theorem~\ref{theo:exist} and then Theorem~\ref{theo:SOLA}. 

\begin{proof}[Proof of Theorem~\ref{theo:exist}] The proof is provided in steps based  mainly on ideas of \cite{BCDM,CiMa,DHM}.\newline

\noindent\textbf{Step 1. }\emph{Approximate problems}.  We consider
\begin{equation}
    \label{eq:main:fk}
    -\tedv\opA(x,\teDu_j)=\tef_j
\end{equation}
with
\[\tef_j(x):= \int_\Rm \vr_j(x-y)\, d\temu(y),\]
where $\vr_j$ stands for a standard mollifier i.e. for a nonnegative, smooth, and even function such that $\int_\r\vr(s){\rm\,d}s=1$ we define  $\vr_j(s)=j\vr(js)$ for $j\in \N$. Of course   \begin{equation}
    \label{fj-conv-mu}
\tef_j\xrightharpoonup[]{*}\temu\qquad\text{and}\qquad\sup_j\|\tef_j\|_{L^1(\Omega)}\leq|\temu|(\rn)<\infty.
\end{equation}
By Remark~\ref{rem:weak-sol} 
 one finds $\teu_j\in W_0^{1,G}(\Omega,\Rm)$
such that for every $\tevp\in W_0^{1,G}(\Omega,\Rm)$ it holds that
\begin{equation}
    \label{408bis}
\int_\Omega \opA(x,\teDu_j): \teDvp\,dx=\int_\Omega\tevp\,\tef_j dx.
\end{equation} 
 We need to show that there exists $\teu\in  {\Upsilon}^{1,G}_0(\Omega)$ and such that (up to a subsequence) there holds a convergence  \begin{equation}
     \teu_j\to \teu \qquad \text{a.e. in $\Omega$.}\label{conv:ujtou-a.e.}
 \end{equation} 

\noindent\textbf{Step 2.  }\emph{The first a priori estimate}. Let us recall~\eqref{DTk} and structure assumption~\eqref{ass:struct} imposed on $\opA$. Making use of \eqref{ass:gr'n'coer},  $\tevp=T_t(\teu_j)$ in \eqref{408bis} and by~\eqref{fj-conv-mu}  for $t>0$ 
and for any $j\in \N$ one obtains that
\begin{align}\label{414}
c \int_{\Omega} G&(|D \Theta_{t/2}(\teu_j)|)\, dx
\leq 
\int_{\{|\teu_j|<t\}}G(|\teDu_j|)\, dx \\
&\nonumber \le \frac 1{c_1^\opA }\int_{\{|\teu_j|<t\}}\opA(x, \teDu_j):\teDu_j\, dx \le \frac{1}{c_1^\opA } 
 \int_{\Omega}\opA(x, \teDu_j):DT_t(\teu_j)\, dx 
\\ \nonumber  & 
=\frac{1}{c_1^\opA }
\int_{\Omega} \tef_j \, T_t(\teu_j)\, dx\le \frac{|\temu|(\Omega)}{c_1^\opA }t=Mt\,,
\end{align}
Here, $\Theta_t$ is any $C^1_0(\Rm,\Rm)$ function satisfying \eqref{Theta-t}. Due to this choice $c>0$ is independent of $t$.

\noindent{\bf Step 3.  }\emph{Almost everywhere convergence of $\{\teu_j\}$}. 
Now we will show that for fixed $\ve>0$, there exists $t_1=t_1(t_0,\ve,\data)$ (and not on $j$), such that 
\begin{equation}
    \label{uj-vanishes}
\big|\{|\teu_j|\geq t\}\big| <\ve\quad\text{ for } t \geq t_1.
\end{equation}
We prove it in two cases -- when $G$ is growing at infinity quickly and slowly in the sense of \eqref{intG}. If necessary for \eqref{int0G}, we substitute $G$ with $G^0$ defined in~\eqref{modification-at-0} and assume that $t_1\geq 1$.
 
  If \eqref{intG}$_1$ holds, then we will show that there exists a  constant $c=c(n,\sigma,|\Omega|)$  such that for  {$\sigma \geq n$} 
\begin{equation}\label{4.4}
|\{|\teu_j|\geq t\}|\leq \frac{Mt}{G_\sigma(ct^{\frac 1{\sigma'}}/M^{\frac 1\sigma})} \qquad \hbox{for $t >0$.}
\end{equation}
Clearly, $|\teu_j| \in W^{1,G}(\Omega)$, and hence $\T_t(|\teu_j|) \in W^{1,G}(\Omega)$ for $t>0$.  
 By  the Orlicz--Sobolev inequality \eqref{Sob-slow} applied  to the function $\T_t (|\teu_j|)$,
\begin{equation}\label{os1}
\int _\Omega G_\sigma\Bigg(\frac{\T_t (|\teu_j|)}{C \big(\int _\Omega G(|D\T_t (|\teu_j|)|)dy\big)^{1/\sigma}}\Bigg)\, dx \leq \int _\Omega G(|D\T_t (|\teu_j|)|)dx.
\end{equation}
One has that
\begin{equation}\label{os4}
\int _\Omega G(|D \T_t (|\teu_j|)|)dx \leq \int _{\{|\teu_j|< t\}} G(|D\teu_j|)dx \quad \hbox{for $t>0$,}
\end{equation}
and $ \{|\T_t(|\teu_j|)|\geq t\} =  \{|\T_t(|\teu_j|)|= t\} = \{|\teu_j|\geq t\} $ for $t>0$. Thus, 
\begin{align}\nonumber
|\{|\teu_j|\geq t\}| G_\sigma&\bigg(\frac{ t}{C (\int _{\{|\teu_j|< t\}} G(|D\teu_j|)dy)^{1/\sigma} }\bigg )\\
\nonumber& \leq \int _{\{|\teu_j|\geq t\}} G_\sigma\Bigg(\frac{\T_t(|\teu_j|)}{C \big(\int _{\{|\teu_j|< t\}} G(|D \teu_j|)dy\big)^{1/\sigma}}\Bigg)\, dx \\ \label{os2}& \leq \int _{\{|\teu_j|< t\}} G(|D \teu_j|)dx 
\end{align}
for $t>0$.  Hence, by \eqref{414},
\begin{align}\label{os3}
|\{|\teu_j|\geq t\}| G_\sigma\bigg(\frac{ t}{C (M t)^{1/\sigma}}\bigg)  \leq M t \qquad \hbox{for $t > 0$.}
\end{align}
Hence~\eqref{4.4} follows in the case \eqref{intG}$_1$.
Since $t\mapsto G(t)/t^{i_G}$ is nondecreasing and $i_G>1$, it holds $\lim _{t \to \infty} {G_\sigma(t)}/{t^{\sigma'}} = \infty.$
Note that~\eqref{uj-vanishes} is a direct consequence of~\eqref{4.4}.
 
If    \eqref{intG}$_2$ holds, then
 there exists  a constant  $t_1=t_1(t_0, \sigma, M)$ such that
\begin{equation}\label{4.4inf}
|\{|\teu_j|\geq t\}|=0  \qquad \hbox{for $t \geq t_1$.}
\end{equation} By Orlicz--Sobolev embedding~\eqref{Sob-fast} we have
\begin{align*}
\|\T_t (|\teu_j|)\|_{L^\infty(\Omega)} & \leq c \,J_n\bigg(\int _\Omega G(|D\T_t (|\teu_j|)|)\, dx\bigg)
\\ \nonumber & = c \,J_n\bigg(\int _{\{|\teu_j|< t\}} G(|\teDu_j|)\, dx\bigg) \leq c \, J_n(Mt) \qquad \hbox{for $\ t > 0$.}
\end{align*}
Observe that $J_n$ is concave, therefore when we set $t_1$ such that $t_1=cJ_n(Mt_1)=\frac{F_n(cM)}{M}$ we have 
\begin{align}\label{os8}
|\{|\teu_j|>t\}|  =0 \qquad \hbox{if $\ t > t_1$}.
\end{align}

\medskip

{We will show that $\{\teu_j\}$ is a Cauchy sequence in measure and there exists a measurable function 
$\teu: \Omega \to \Rm$ such that, up to subsequences, we have}
\begin{equation}\label{408ter}
\teu_j\rightarrow \teu \qquad \hbox{a.e. in }\Omega.
\end{equation}
 
Fix any coordinate $i\in\{1,\dots, m\}.$ Given any $t, \tau>0$, one has that
\begin{equation}\label{409}
|\{|\teu_k - \teu_m|>\tau \}|\le |\{| \teu_k |>t \}|+|\{| \teu_m |>t \}|+ |\{| \Theta_t(\teu_k) - \Theta_t(\teu_m)|>\tau  \}|,
\end{equation}
for $k,m\in \N$. Having~\eqref{uj-vanishes}, for any fixed $\varepsilon>0$,
the number $t$ can be chosen so large that
\begin{equation}\label{412}
 |\{| \teu_k|>t \}|<\varepsilon  \quad \hbox{and} \quad  |\{| \teu_m|>t \}|<\varepsilon
\end{equation} 
for every $k, m \in \N$.
By Orlicz--Sobolev inequality \eqref{Sob-slow}-\eqref{Sob-fast} and a priori estimate \eqref{414}, the sequence $\{\Theta_t(\teu_j)\}_j$ is
bounded in $W^{1,G}(\Omega, \Rm)\subset W^{1,1}(\Omega,\Rm)$. By the compactness of embedding $W^{1,1}(\Omega, \Rm)\Subset L^1(\Omega, \Rm)$,  the sequence $\{\Theta_t(\teu_j)\}$
converges (up to subsequences) to some function in
$L^1(\Omega)$. In particular, $\{\Theta_t(\teu_j)\}$ is a Cauchy
sequence in measure in $\Omega$. Thus,
\begin{equation}\label{416}
|\{ |\Theta_t(\teu_k)-\Theta_t(\teu_m)|>\tau  \}|\le  \varepsilon
\end{equation}
provided that $k$ and $m$ are sufficiently large. By \eqref{409},
\eqref{412} and \eqref{416}, $\{\teu_j\}$ is (up to a subsequence) a
Cauchy sequence in measure in $\Omega$, and hence there exists a
measurable function $\teu: \Omega\rightarrow \R$ such that \eqref{408ter}
holds. \newline

\noindent\textbf{Step 4.  }\emph{A gradient estimate}.  We will use the ideas of~\cite{BCDM} to show that
\begin{equation}\label{andrea12}
\int_\Omega g (|\teDu_j|)\,dx\leq c |\temu|(\Omega)\qquad 
\text{for some constant $c=c(n,|\Omega|)>0$.} 
\end{equation}
In the case of \eqref{intG}$_2$ for $n=\sigma$, Sobolev embedding \eqref{Sob-fast}  directly gives $F_{n}^{-1}\left(\int_{\Omega} G(|D\teu_j|) \, dx \right) \leq c |\temu|(\Omega)$, so from now on we concentrate on the slowly growing $G$ (i.e. satisfying \eqref{intG}$_1$).

From inequality \eqref{414} we deduce that
\begin{align}\nonumber
|\{G(|\teDu_j|)>s, |\teu|\leq t\}| &\leq     \frac 1s \int_{\{G(|\teDu_j|)>s, |\teu_j|\leq t\}} G(|\teDu_j|)\, dx \\
&\leq \frac{t \|\tef_j\|_{L^1(\Omega, \Rm)}}s\leq \frac{t |\temu|(\Omega)}s  \quad \hbox{for $t >0$ and $s > 0$.}\label{lem2.30}
\end{align}
Let $\sigma\geq n$. 
By combining \eqref{lem2.30} and \eqref{os3} we get that for $t >0$ and $s>0$ it holds
\begin{align}\label{lem2.4}
|\{G(|\teDu_j|)>s\}| & \leq  |\{|\teu_j| >t\}| + |\{G(|\teDu_j|)>s,|\teu_j|\leq t\}| \\ \nonumber & \leq \frac{t \|\tef_j\|_{L^1(\Omega, \Rm)}}{G_\sigma(Ct^{\frac 1{\sigma'}}/{ \|\tef_j\|_{L^1(\Omega, \Rm)}^{\frac 1\sigma})}} + \frac{ t \|\tef_j\|_{L^1(\Omega, \Rm)}}s.
\end{align}
When we recall that $\|\tef_j\|_{L^1(\Omega, \Rm)}\leq|\temu|(\Omega)$, the choice $t = \big(\tfrac 1C |\temu|(\Omega)^{1/\sigma} G_\sigma^{-1}(s)\big)^{\sigma '}$ in inequality  \eqref{lem2.4} results in
\begin{equation}\label{lem2.5}
|\{{G}(|\teDu_j|)>s\}| \leq  \  \frac{2  |\temu|(\Omega)^{\sigma'}}{C^{\sigma'} }\frac{G_\sigma^{-1}(s)^{\sigma '}}{s} \quad \hbox{for $s>0$.}
\end{equation}
Next, 
set   $s= G(g^{-1}(\tau))$ in \eqref{lem2.5}  and make use of \eqref{GN} to obtain that
\begin{equation}\label{lem2.6}
|\{g(|\teDu_j|)>\tau\}| \leq  \frac{2  |\temu|(\Omega)^{\sigma'}}{C^{\sigma'} } \frac{H_\sigma(g^{-1}(\tau ))^{\sigma '}}{G(g^{-1}(\tau ))} 
\qquad \hbox{for $\tau >0$,}
\end{equation}
where $H_\sigma$ is defined as in \eqref{GN}.
Thanks to inequality \eqref{lem2.6}, 
\begin{align}\nonumber
\int_\Omega 
 g(|\teDu_j|) \, dx & = \int _0^\infty |\{g(|\teDu_j|)>\tau\}|\, d\tau 
\\ &\leq  \lambda |\Omega| +   2  C^{-\sigma'} |\temu|(\Omega)^{\sigma'}\int_\lambda ^\infty   \frac{H_\sigma(g^{-1}(\tau ))^{\sigma '}}{G(g^{-1}(\tau ))}\, d\tau \label{may99}
\end{align}
for $\lambda >0$. Owing to Lemma~\ref{lem:equivalences},  and to Fubini's theorem, the following chain holds:
\begin{align}\nonumber
\int_\lambda ^\infty   &\frac{H_\sigma(g^{-1}(\tau ))^{\sigma '}}{G(g^{-1}(\tau ))}\, d\tau  \leq c
\int_{g^{-1}(\lambda)} ^\infty   \frac{H_\sigma(s)^{\sigma '}}{sG(s)} g(s) ds\\
&\nonumber= c
\int_{g^{-1}(\lambda)} ^\infty   \frac{g(s)}{sG(s)}  \int _0^s \bigg(\frac t{G(t)}\bigg)^{\frac 1{\sigma-1}}\,dt \,ds
\\ \nonumber & \leq c
\int_{g^{-1}(\lambda)} ^\infty   \frac{1}{s^2}  \int _0^s \bigg(\frac t{G(t)}\bigg)^{\frac 1{\sigma-1}}\,dt \,ds
\\ \nonumber &= c\bigg(\int_0^{g^{-1}(\lambda)}  \bigg(\frac t{G(t)}\bigg)^{\frac 1{\sigma-1}} \int_{g^{-1}(\lambda)} ^\infty   \frac{ds}{s^2} \, dt + \int _{g^{-1}(\lambda)} ^\infty \bigg(\frac t{G(t)}\bigg)^{\frac 1{\sigma-1}} \int_{t} ^\infty   \frac{ds}{s^2} \, dt\bigg)
\\ \nonumber &= c\bigg(\frac{1}{g^{-1}(\lambda)}   \int_0^{g^{-1}(\lambda)}  \bigg(\frac t{G(t)}\bigg)^{\frac 1{\sigma-1}} \, dt + \int _{g^{-1}(\lambda)} ^\infty \bigg(\frac t{G(t)}\bigg)^{\frac 1{\sigma-1}}   \frac{dt}{t}\bigg)
\\\label{may102} &\leq c\bigg(\frac{1}{g^{-1}(\lambda)}   \int_0^{g^{-1}(\lambda)}  \bigg(\frac 1{g(t)}\bigg)^{\frac 1{\sigma-1}} \, dt + \int _{g^{-1}(\lambda)} ^\infty \bigg(\frac 1{g(t)}\bigg)^{\frac 1{\sigma-1}}   \frac{dt}{t}\bigg)=:c(I_1+I_2)
\end{align}
for $c=c(s_G)>0$ and $\lambda >0$. We restrict ourselves to $\sigma>\max\{s_G+1,n\}$. The function $\frac {t^{s_G+\ve}}{g(t)}$
is increasing for every $\ve>0$. Hence, if $0<\ve < \sigma - s_G -1$, then   for $\lambda > 0$ we estimate
\begin{align}\nonumber
I_1&=\frac{1}{g^{-1}(\lambda)}  \int_0^{g^{-1}(\lambda)}  \bigg(\frac 1{g(t)}\bigg)^{\frac 1{\sigma-1}} \, dt  =
\frac{1}{g^{-1}(\lambda)}   \int_0^{g^{-1}(\lambda)}  \bigg(\frac {t^{s_G+\ve}}{g(t)}\bigg)^{\frac 1{\sigma-1}} t^{- \frac{s_G+\ve}{\sigma -1}}\, dt \\ 
\label{may103} & \leq \frac{1}{g^{-1}(\lambda)}  \bigg(\frac {g^{-1}(\lambda)^{s_G+\ve}}{\lambda}\bigg)^{\frac 1{\sigma-1}}   \int_0^{g^{-1}(\lambda)}  t^{- \frac{s_G+\ve}{\sigma -1}}\, dt = \frac{\sigma -1}{\sigma - s_G-1-\ve} \lambda ^{-\frac 1{\sigma -1}}.
\end{align}
On the other hand, if $0<\ve < i_G$, then the  function 
$\frac {t^{\ve}}{g(t)}$
is decreasing. Hence, for $\lambda > 0$
\begin{align}\nonumber
I_2&= \int _{g^{-1}(\lambda)} ^\infty \bigg(\frac 1{g(t)}\bigg)^{\frac 1{\sigma-1}}   \frac{dt}{t}  
 = \int _{g^{-1}(\lambda)} ^\infty \bigg(\frac {t^\ve}{g(t)}\bigg)^{\frac 1{\sigma-1}}  t^{-\frac \ve{\sigma -1}-1}\, dt\\
\label{may104} &\leq  \bigg(\frac {g^{-1}(\lambda)^\ve}{\lambda}\bigg)^{\frac 1{\sigma-1}}  \int _{g^{-1}(\lambda)} ^\infty t^{-\frac \ve{\sigma -1}-1}\, dt
 = \frac{\sigma -1}\ve \lambda ^{-\frac 1{\sigma -1}}.
\end{align}
Inequalities \eqref{may102}--\eqref{may104} entail that there exists a constant $c=c(\sigma, i_G, s_G)$ such that for $\lambda >0$
\begin{align}\label{may105}
\int_\lambda ^\infty   \frac{H_\sigma(g^{-1}(\tau ))^{\sigma '}}{G(g^{-1}(\tau ))}\, d\tau \leq c \lambda ^{-\frac 1{\sigma -1}}.
\end{align}
Inequality \eqref{andrea12} follows from \eqref{may99} and \eqref{may105}, with the choice $\lambda = |\temu|(\Omega)$. \newline

\noindent{\bf Step 5. } \emph{Approximate differentiability of $\teu$}. Let us fix $t>0$ and the function $\Theta_t$.
By \eqref{414}, the sequence $\{ \Theta_t(\teu_j)\}_j$ is bounded in $W^{1,G}(\Omega;\R^m)$. Since, by \eqref{408ter},  we already know that there exists a measurable function $\teu$ such that
$$
\teu_j \to \teu \qquad \text{ a.e.} 
$$
then, after passing to a subsequence we may assume that
$$
\Theta_t(\teu_j) \xrightharpoonup[j \to \infty]{} \Theta_t(\teu) \qquad \text{in $W^{1,G}(\Omega;\Rm)$}.
$$
In particular, $\Theta_t(\teu)$ is a Sobolev function in $\Omega$, and thus it is approximately differentiable a.e. in $\Omega$ (see e.g. \cite{EG} or \cite{Ziem}). As $\teu$ is measurable, it is also approximately continuous a.e. in $\Omega$. 
{Recall that $
M_{t/2}= \{ x \in \Omega : |\teu(x)| <t/2 \}.
$}
Let $x_0 \in M_{t/2}$ be a point of approximate continuity of $\teu$ and approximate differentiability of $\Theta_t(\teu)$. It follows that
$$
\lim_{r\to 0}\frac{\left|\left\{x\in B(x_0,r):\ |\teu(x)-\teu(x_0)|\geq t/4\right\}\right|}{r^n}=0
$$
and for every $\ve >0$
$$
\lim_{r\to 0}
\frac{\left|\left\{x\in B(x_0,r):\ |\Theta_t(\teu(x))-\Theta_t(\teu(x_0))-D(\Theta_t\circ \teu)(x_0)(x-x_0)|>\ve r\right\}\right|}{r^n}=0\,.
$$
Consider a set
$$
E_{r,\ve} = \{ x \in B(x_0,r) : | \teu(x) - \teu(x_0) - D(\Theta_t \circ \teu)(x_0)(x-x_0) | > \ve r \}.
$$
Approximate continuity of $\teu$ implies
\begin{equation*}
\lim_{r\to 0}\frac{|E_{r,\ve} \cap \{| \teu(x)-\teu(x_0)|\geq t/4 \}|}{r^n}=0\,.
\end{equation*}
On the other hand, if  {$x_0 \in M_{t/2}$} and $| \teu(x)-\teu(x_0)|<t/4$ then $|\teu(x)| < 3/4t$.
Since $\Theta_t \equiv \id$ on $B(0,t)$ it follows that if $| \teu(x)-\teu(x_0)|<t/4$ we have $\Theta_t(\teu(x)) = \teu(x)$, and by approximate differentiability of $\Theta_t(\teu)$ at $x_0$ we have
\begin{equation*}
\lim_{r\to 0}\frac{|E_{r,\ve} \cap \{| \teu(x)-\teu(x_0)|<t/4 \}|}{r^n}=0\,.
\end{equation*}
Hence
$$
\lim_{r\to 0}\frac{|E_{r,\ve}|}{r^n}=0\,.
$$
which means that $\teu$ is approximately differentiable at $x_0$, and
$$
{\rm ap} D\teu(x_0) = {\rm ap} D(\Theta_t \circ \teu)(x_0) = D(\Theta_t \circ \teu)(x_0). 
$$
It follows that $\teu$ is approximately differentiable for almost every $x \in \Omega$ because $ \bigcup_{t>0} M_t = \Omega \setminus E$ where $E$ is a set of measure zero. Moreover,  ${\rm ap}\teDu(x) = \tens{Z_u}(x)$ where $\tens{Z_u}$ is the generalized gradient defined in \eqref{gengrad'} for a.e. $x \in \{ |\teu| \leq t\}$.\newline

\noindent{\bf Step 6.  }\emph{Limits in measure}.  We will show that the limit function $\teu$ satisfies $\teu\in \Upsilon^{1,G}(\Omega)$ {and} up to a subsequence
\begin{equation}\label{425}
\teDu_j \rightarrow \teDu \qquad \hbox{ in measure}\,
\end{equation}
and 
\begin{equation}\label{425'}
\opA(x,\teDu_j) \rightarrow \opA(x,\teDu) \qquad \hbox{ in measure}\,,
\end{equation}
where $\teDu$ is the generalized
gradient of $\teu$ in the sense of the function $\tens{Z_u}$ appearing in \eqref{gengrad'}.

We apply the results of~\cite[Section~5]{DHM-gen} recalled and commented in Appendix. There exists a sequence $\{\alpha_k\}$ diverging to infinity, such that the limit function $\teu\in W^{1,1}(M_{\alpha_k})$ for each $k$ if $M_{\alpha_k}:=\{x\in\Omega:\ |\teu|<\alpha_k\}$.  It is enough for convergence $\teDu_j\to \teDu$ a.e. to prove it for each $k$. Fix $\alpha=\alpha_k$ for arbitrary $k$. By~\eqref{414} on $M_\alpha$ function $\teu$ is weakly differentiable. Due to \cite[Chapter 6]{Ziem} or \cite[Chapter 3]{EG} on $M_\alpha$  the weak gradient of $\teu$ satisfies also the definition of the approximate gradient ${\rm ap}\teDu$. In order to prove that the limit is indeed $\teDu$, we employ Lemmas~\ref{Lem11} and~\ref{Lem12}. Their consequence is that $\opA(\cdot,\teDu_j)\to\opA(\cdot,\teDu)$ strongly in  $L^1(M_\alpha,\Rm)$ and $\teDu_j\to{\rm ap}\teDu=\teDu$ in measure in $\Omega$. We can employ Lemma~\ref{Lem11} as the structural assumptions {\it (A1)--(A6)} are satisfied. Indeed, we assume that $\opA$ is a Carath\'eodory's function, which is monotone by~\eqref{ass:str-mon} and such that~\eqref{ass:struct}. Moreover, for each $j$ we know that $\teu_j$ is weakly differentiable in $M_\alpha$, $\tef_j=-\dv\opA(\cdot,\teDu_j)$ by~\eqref{eq:main:fk}, $\{\tef_j\}_j$ is bounded in $L^1$ by~\eqref{fj-conv-mu}, by Step~2 we know that $\teu_j\to\teu$ in measure, and -- as already mentioned -- the limit function $\teu$ has an approximate gradient ${\rm ap}\teDu$. Furthermore, $\{\opA(\cdot,\teDu_j)\}_j$ is uniformly integrable and $\{\opA(\cdot,\teDu_j)\}_j$ is uniformly bounded in $L^1$ by~\eqref{ass:gr'n'coer}  and~\eqref{andrea12}.  The final argument identifying the limit is provided by Lemma~\ref{Lem12} due to the strong monotonicity of the operator~\eqref{ass:str-mon}. Since ${\rm ap}\teDu=\tens{Z_u}$ in $M_\alpha$ and $\alpha=\alpha_k$ for an arbitrary $k$, we can conclude with~\eqref{425} and \eqref{425'}. \newline

\color{black}

\par\noindent
{\bf Step 7.}  \emph{Approximable solution $\teu$ shares regularity of claims {\rm (i)} and {\rm (ii)}.}

\par  In {\rm (ii)} $G$ grows so fast that \eqref{intG}$_2$ holds. By~\eqref{Sob-fast}  $W^{1,G}(\Omega,\Rm)\subset L^\infty(\Omega,\Rm)$. Since we know that~\eqref{distr} holds true, to prove that an approximable solution $\teu$ is a weak solution it is enough to show that it belongs to $W^{1,G}(\Omega,\Rm)$. From~\eqref{4.4inf} we directly deduce that 
\[\|\teu_j\|_{L^\infty(\Omega,\Rm)}\leq t_1\]
with $t_1$ independent of $j$. By a priori estimate~\eqref{414} we have
\[\int_\Omega G(|\teDu_j|)\,dx=\int_{\{|\teu_j|\leq t_1\}} G(|\teDu_j|)\,dx\leq Mt_1.\] Having~\eqref{425}, we get that also \[\int_\Omega G(|\teDu|)\,dx\leq Mt_1\] and {\rm (ii)} follows.

Let us concentrate on {\rm (i)}. We deal with $G$ satisfying \eqref{intG}$_1$. Since $\teu_j$ satisfies~\eqref{4.4} with constants independent of $j$ and we have convergence in measure $\teu_j\to \teu$ from Step~3, it suffices to take $\sigma=n$ to get 
\begin{equation}\label{marc10}
|\{|\teu|>t\}| \vt_n(t/c)\leq c
\end{equation}
for a suitable positive constant $c$ and sufficiently large $t$.  Hence, $|\teu| \in L^{\vt_n, \infty}(\Omega)$. 
\\ Since for $\teDu_j$ inequality \eqref{lem2.5} holds for a suitable $M$ and sufficiently large $s$ with constants independent of $j$ and we have convergence in measure $\teu_j\to \teu$ from Step~5, it suffices to take $\sigma=n$ to get  
$$|\{|\teDu|>t\}| \theta_n (t)\leq c$$
for a suitable constant $c$ and sufficiently large $t$. This implies that $|\teDu| \in L^{\theta_n, \infty}(\Omega)$.  
\end{proof}

\begin{proof}[Proof of Theorem~\ref{theo:SOLA}] We follow ideas of~\cite{Y-new}. SOLAs and approximable solutions can be obtained in the same approximating procedure as in Step~1 of the proof of Theorem~\ref{theo:exist}. Therefore, they share all properties proven in Steps~{2--6} of the same proof.

We concentrate on showing that an approximable solution is a SOLA. For a fixed approximable solution $\teu$ we might construct an approximate sequence $(\teu_j)$ as in Steps 1 of the proof of Theorem~\ref{theo:exist}.   
We start with showing that
\begin{equation}\label{L1-conv} \{|\teDu_j|\} \text{ is a Cauchy sequence in } L^1(\Omega). \end{equation}
Let us recall that $\Psi_n$ is defined in~\eqref{SOLA-Psi}. Writing $t = G(s)$ in \eqref{lem2.5} and taking $\sigma = n$, we have
\[|\{|\teDu_j|>t\}| \leq \frac{c |\temu|(\Omega)^{n'}}{\Psi_n (t)}.\]
Then for any $k,m \in \mathbb{N}$,
\[
|\{|\teDu_k - \teDu_m|>t\}|
\leq |\{|\teDu_k|>t/2\}| + |\{|\teDu_m|>t/2\}|
\leq \frac{c |\temu|(\Omega)^{n'}}{\Psi_n (t/2)}.
\]
For any constants $0< \ve_0 <M_0$ chosen in a moment, 
\begin{align*}
    \int_{\Omega} |\teDu_k - \teDu_m| \, dx & = \int_{0}^{\ve_0} |\{ |\teDu_k - \teDu_m| > t\}| \, dt + \int_{\ve_0}^{M_0} |\{ |\teDu_k - \teDu_m| > t\}| \, dt \\
    & \quad + \int_{M_0}^{\infty} |\{ |\teDu_k - \teDu_m| > t\}| \, dt \\
    & \leq \ve_0 |\Omega| + M_0 |\{ |\teDu_k - \teDu_m| > \ve_0\}| + c |\temu|(\Omega)^{n'} \int_{M_0}^{\infty} \frac{1}{\Psi_n(t)} \, dt.
\end{align*}
For any $\ve>0$, we now take $M_0$ large enough to make the last integral less than $\ve$ due to \eqref{SOLA-Psi}.
Taking $\ve_0=\ve/(1+|\Omega|)$ and then using \eqref{425}, we conclude \eqref{L1-conv}.

Note  {that by~\eqref{425'} $\{ \opA(x,\teDu_j)\}_j$ is a~Cauchy sequence in measure. }  We have
\begin{align*}
    \int_{\Omega}& |\opA(x,\teDu_k) - \opA(x,\teDu_m)| \, dx\\
    & = \int_{0}^{\ve_0} |\{ |\opA(x,\teDu_k) - \opA(x,\teDu_m)| > t\}| \, dt  \\
    & \quad + \int_{\ve_0}^{M_0} |\{ |\opA(x,\teDu_k) - \opA(x,\teDu_m)| > t\}| \, dt \\
    & \quad + \int_{M_0}^{\infty} |\{ |\opA(x,\teDu_k) - \opA(x,\teDu_m)| > t\}| \, dt \\
    & \leq \ve_0 |\Omega| + M_0 |\{ |\opA(x,\teDu_k) - \opA(x,\teDu_m)| > \ve_0\}| \\
    & \quad + \int_{M_0}^{\infty} |\{ g(|\teDu_k|) + b(x)> t/2 \}| \, dx  + \int_{M_0}^{\infty} |\{ g(|\teDu_m|) + b(x)> t/2 \}| \, dx \\
    & \leq \ve_0 |\Omega| + M_0 |\{ |\opA(x,\teDu_k) - \opA(x,\teDu_m)| > \ve_0\}|+ 2 \int_{M_0}^{\infty} |\{ b(x) > t/4 \}| \, dx \\
    & \quad + \int_{M_0}^{\infty} |\{ g(|\teDu_k|) > t/4 \}| \, dx + \int_{M_0}^{\infty} |\{ g(|\teDu_m|) > t/4 \}| \, dx \\
    & \leq \ve_0 |\Omega| + M_0 |\{ |\opA(x,\teDu_k) - \opA(x,\teDu_m)| > \ve_0\}|+ \frac{8}{M_0} \int_{\{ b(x) > M_0/4 \}} b(x) \, dx \\
    & \quad + \int_{M_0}^{\infty} |\{ g(|\teDu_k|) > t/4 \}| \, dx + \int_{M_0}^{\infty} |\{ g(|\teDu_m|) > t/4 \}| \, dx .
\end{align*}
Using \eqref{lem2.6}, we can make the last line on the right-most-side to be as small as we want by taking $t$  large enough. Moreover, the last integral in the second last line converges to 0 as $M_0 \to \infty$, so this term can also be made arbitrarily small. We may take $\ve_0=\ve/|\Omega|$ and then use the fact that $\{\opA(x,\teDu_j)\}_j$ is a Cauchy sequence in measure to make the remaining term small. In the end, we get that  $\{\opA(x,\teDu_j)\}_j$ is a Cauchy sequence in the sense of $L^1$-convergence  {having a limit as in~\eqref{425'}}.
\end{proof}

{
\section{Appendix}

In this appendix, we recall the fundamental theorem about Young measures and lemmas in \cite{DHM-gen} that we have used in the proof of Theorem \ref{theo:exist}. Let us mention that the notion of Young measures dates back to~\cite{Young} whereas a relevant reference for the proofs of their properties is~\cite{Ball-Young}. 
 
Let  $C_0^0 (\Rm) := \{\vp\in C^0(\Rm):\ \lim_{|z|\to\infty}|\vp(z)| = 0\}$.

\begin{theo}[Fundamental theorem for Young measures]\label{theo:Youngmeas} 
Let $\Omega \subset\rn$ and $\texi_j : \Omega \to\Rm$ be a~sequence of measurable functions. Then there exists a~subsequence $\{\texi_{j,k}\}$ and a~family $\{\nu_x\}_{x\in\Omega}$ of non-negative Radon measures on $\Rm$, such that:
\begin{itemize}
\item[i) ]$\|\nu_x\|_{\mathcal{M}(\Rm)}=\int_{\Rm} \, d\nu_x \leq 1$ for a.e. $x \in \Omega$.
\item[ii) ] For $\temf \in C_0^0(\Rm)$, we have $\temf (\texi_{j,k})\xrightharpoonup[]{*} \bar{\temf}$ weakly-$\ast$ in $L^\infty(\Omega)$ and  $\bar{\temf}(z)
=\langle \nu_x, f\rangle$.
\item[iii) ] If for all $R>0$ it holds $\lim_{L\to\infty}\sup_{k\in\N}|\{x\in\Omega\cap B(0,R):\ |\texi_{k}(z)|\geq L\}|=0$, then $\|\nu_x\|_{\mathcal{M}(\Rm)}=1$ for a.e. $x \in \Omega$, and for all measurable $A\subset \Omega$ there holds $\temf(\texi_{j,k})\xrightharpoonup[]{} \bar{\temf}=\langle \nu_x, f\rangle$ in $ L^1(A)$ for continuous $f$ provided $\{f(\texi_{k})\}_k$ is weakly precompact in $L^1(A)$.
\end{itemize}  
The family of maps $\nu_x : \Omega \to \mathcal{M}(\Rm)$ is called the
Young measure generated by~$\{\texi_{k}\}_k$.
\end{theo}

By iii) we infer that if $A\subset\Omega$ and $f:A\times\Rm\to \R$ is a Carath\'eodory's function, then $f(\cdot,\texi_k)\xrightharpoonup{}\langle \nu_x, f(x,\cdot)\rangle$ provided the sequence $\{f(\cdot,\texi_{k})\}_k$ is weakly precompact in $L^1(A)$. Moreover, if $|\Omega|<\infty$, then \[\text{$\texi_k\to\texi$ in measure
$\iff$ the Young measure assosciated to $\texi_k$ is $\delta_{\texi(x)}$.}\]
If $\nu=\{\nu_x\}_{x \in \Omega}$ is a Young measure associated with a sequence of gradients $\{D \texi_k\}$, where $\{\texi_k\}$ is bounded in $W^{1,G}(\Omega)$ for any $N$-function $G$, then we call $\nu$ a $W^{1,G}$-gradient Young measure. If $\{\nu_x\}_{x \in \Omega}$ is a $W ^{1,G}$-gradient Young measure, then there exists a function
$\teu \in W^{1,G}(\Omega;\Rm)$, such that $\teDu(x ) = \langle\nu_x , \mathsf{Id}\rangle$ almost everywhere.

By the arguments of \cite[Lemma~8 and~9]{DHM-gen}, we infer the following fact.
\begin{lem}\label{Lem8-9} Suppose $G\in\Delta_2\cap \nabla_2$.
    If $\{\teu_k\}_k\in W^{1,1}(\Omega;\Rm)$, for any $\alpha>0$ it holds $\sup_k\int_{|\teu_k|\leq\alpha\}} G(|\teDu_k|)\,dx\leq C(\alpha)<\infty$ and for $s>0$ it holds $\sup_k\int_\Omega |\teu_k|^s\,dx\leq C<\infty,$ then
    up to a subsequence $\teu_k\to \teu$ in measure and the Young measure  $\nu_x$ generated by a subsequence of $\{\teDu_k\}_k$ satisfies  ${\rm ap}\teDu= \langle \nu_x, \mathsf{Id} \rangle$.
\end{lem}

In the next lemmas, we consider $\opA,\teu_k,$ and $\teu$ satisfying the following assumptions.
\begin{enumerate}
\item[\it (A1)] $\opA : \Omega \times \rnm \to \rnm$ is a Carath\'eodory function satisfying \eqref{ass:struct}.
\item[\it (A2)] $\teu_k \in W^{1,1}(\Omega,\Rm)$ and $\int_{\Omega} |D \teu_k|^{s} \, dx \leq C$ uniformly in $k$ for some $s>0.$
\item[\it (A3)] The sequence $\opA_{k}(x) = \opA(x, D \teu_k(x))$ is uniformly integrable.
\item[\it (A4)] The sequence $\teu_k$ converges in measure to $\teu,$ and $\teu $ is approximately differentiable a.e. in $\Omega.$
\item[\it (A5)] The sequence $f_k := - \dv \opA_k$ is bounded in $L^1(\Omega)$.
\item[\it (A6)] $D \teu_k \in L^G_{loc}$ and $\opA_k \in L^{\wt{G}}_{loc}$ for some $N$-function $G\in\Delta_2\cap\nabla_2$.
\end{enumerate}

\begin{lem}\label{Lem11}Let (A1)-(A6) be in force. Then $\opA_k$ converges (up to subsequences) weakly in $L^1(\Omega)$ to $\bopA$, which is given by $\bopA(x) = \langle \nu_x, \opA(x, \cdot) \rangle$.
Moreover, we have the inequality
$$\int_{\rnm} \opA(x,\lambda) : \lambda \, d \nu_x (\lambda) \leq \bopA(x) : \apdu(x) \quad \text{for a.e. } x \in \Omega. $$
\end{lem}
\begin{proof}[Comments on the proof]
    The proof of the above lemma is the same as \cite[Lemma~11]{DHM-gen} under assumptions \cite[(5.1)-(5.7)]{DHM-gen}. It only remains to explain that \cite[(5.7)]{DHM-gen} reading  $\teDu_j\in L^r$ and $\opA(x,\teDu_j)\in L^{r'}$ for some $r\in (1,\infty)$ can be changed to~{\it (A6)}. Note that \cite[(5.7)]{DHM-gen} was not directly used in the proof of \cite[Lemma~11]{DHM-gen}, and it can be easily passed by in our situation. Indeed, it is only applied to make $\opA(x,\teDu_j):\teDu_j\in L^1$ and between (5.13) and (5.14) to justify passing to the limit of a product $f_j\,g_j\to fg$ knowing $f_j\xrightharpoonup{} f$ in $L^r$, $r\geq 1$ and $g_j,g\in L^\infty$, $g_j\to g$ almost everywhere. The first property holds due to {\it (A6)}. 
    The second one holds for $r=1$ by {\it (A5)}. 
    The rest of the proof does not need any modification.
\end{proof}

As a consequence of Lemma~\ref{Lem11}, the next one follows precisely the same lines as \cite[Lemma~12]{DHM-gen}. 
\begin{lem}
   \label{Lem12}Let (A1)-(A6) be in force. Suppose  that assumptions of Lemma~\ref{Lem8-9} hold true and that
 $\opA$ satisfies~\eqref{ass:str-mon}. Then for $\alpha>0$ it holds
\begin{equation*}
    \opA(x,\teDu_k)\to\opA(x,{\rm ap}\teDu)\ \text{ strongly in } L^1(\{|\teu|<\alpha\})\ \text{and $\ \teDu_k\to{\rm ap}\teDu$ in measure}\,.
\end{equation*}
\end{lem}
\begin{proof} By Lemma~\ref{Lem8-9} we known that      ${\rm ap}\teDu= \langle \nu_x, \mathsf{Id} \rangle$. Note that having {\it (A4)} we can infer that {\it (A3)} is satisfied only on $\{|\teu|<\alpha\}$ for any $\alpha$. Using Lemma~\ref{Lem11} we infer that $\supp \nu_x\subset \{\lambda:\ \bopA( x):\lambda=0\}$. Then the strict monotonicity of the operator $\opA$ implies that $\nu_x=\delta_{{\rm ap}\teDu}$. Consequently, $\bopA(x)=\opA(x,{\rm ap}\teDu(x))$, $\teDu_k\to{\rm ap}\teDu$ in measure, and $\opA(x,\teDu_k)\to\opA(x,{\rm ap}\teDu)$ in $L^1(\{|\teu|<\alpha\})$.
\end{proof}

}
\section*{Acknowledgements} I. Chlebicka is supported by NCN grant no. 2019/34/E/ST1/00120. A. Zatorska-Goldstein is supported by NCN grant no. 2019/33/B/ST1/00535. Y. Youn is supported by NRF grant no. 2020R1C1C1A01009760.

\end{document}